
\parindent=0pt
\parskip=4pt
\input epsf

\raggedbottom
\hsize=5.5truein

%
\chardef\active 13
\chardef\other = 12
\def\deactivate{%
\catcode`\\ =\other
\catcode`\} =\other
\catcode`\{ =\other
\catcode`\& =\other
\catcode`\$ =\other
\catcode`\# =\other
\catcode`\~ =\other
\catcode`\_ =\other
\catcode`\^ =\other
\catcode`\% =\other
}
\def\makeactive#1{\catcode`#1=\active\ignorespaces}
{\makeactive\^^M%
\makeactive\^^K%
\gdef\obeywhitespace{%
\makeactive\^^M%
\makeactive\^^K%
\let^^M = \newline%
\let^^K =\onespace%
\aftergroup\removebox 
\obeyspaces }}
\newdimen\spc   \spc=1.22ex
\newdimen\tbspc \tbspc=4.88ex
\def\newline{\par\indent}

\def\onespace{\hskip1\spc}
\newdimen\savedparskip
\def\removebox{\setbox0=\lastbox}

\def\ttverbatim{\begingroup\tt\deactivate\obeywhitespace}
\def\endverbatim{\endgroup\parskip=\savedparskip{}}


\catcode`\@=\active
\newskip\ttglue \ttglue=.5em plus .2em minus .1em
{\gdef@{\savedparskip=\parskip\parskip=0pt\ttverbatim\spaceskip=\ttglue\let@=\endverbatim}}

\input ./inputs/font.defs
\def\subdef#1{\gdef\globalColor##1{##1}}
%
%
%

%

%
%
\def\newColor #1 {\expandafter\def\csname #1\endcsname##1{##1}%
   \expandafter\def\csname text#1\endcsname{\subdef{#1}%
    }}%
%
%
\newColor GreenYellow     
\newColor Yellow	  
\newColor Goldenrod	  
\newColor Dandelion	  
\newColor Apricot	  
\newColor Peach		  
\newColor Melon		  
\newColor YellowOrange	  
\newColor Orange	  
\newColor BurntOrange	  
\newColor Bittersweet	  
\newColor RedOrange	  
\newColor Mahogany	  
\newColor Maroon	  
\newColor BrickRed	  
\newColor Red		  
\newColor OrangeRed	  
\newColor RubineRed	  
\newColor WildStrawberry  
\newColor Salmon	  
\newColor CarnationPink	  
\newColor Magenta	  
\newColor VioletRed	  
\newColor Rhodamine	  
\newColor Mulberry	  
\newColor RedViolet	  
\newColor Fuchsia	  
\newColor Lavender	  
\newColor Thistle	  
\newColor Orchid	  
\newColor DarkOrchid	  
\newColor Purple	  
\newColor Plum		  
\newColor Violet	  
\newColor RoyalPurple	  
\newColor BlueViolet	  
\newColor Periwinkle	  
\newColor CadetBlue	  
\newColor CornflowerBlue  
\newColor MidnightBlue	  
\newColor NavyBlue	  
\newColor RoyalBlue	  
\newColor Blue		  
\newColor Cerulean	  
\newColor Cyan		  
\newColor ProcessBlue	  
\newColor SkyBlue	  
\newColor Turquoise	  
\newColor TealBlue	  
\newColor Aquamarine	  
\newColor BlueGreen	  
\newColor Emerald	  
\newColor JungleGreen	  
\newColor SeaGreen	  
\newColor Green		  
\newColor ForestGreen	  
\newColor PineGreen	  
\newColor LimeGreen	  
\newColor YellowGreen	  
\newColor SpringGreen	  
\newColor OliveGreen	  
\newColor RawSienna	  
\newColor Sepia		  
\newColor Brown		  
\newColor Tan		  
\newColor Gray		  
\newColor Black		  
\newColor White		  

\subdef{Black}

\input ./inputs/layout.macros
\input ./inputs/satz.macros
\input ./inputs/figure.macros
\input ./inputs/local.macros
\input amssym.def

\def\reflectsto{\;\;\;\longmapsto\;\;\;}

\def\pixel{\vrule width0.75pt height0.375pt depth0.375pt}
\def\dash{\,\raise2.25pt\hbox{\pixel\vrule height0.125pt width12pt depth0.125pt\pixel}\,}
\def\Dash{\,
\raise1.75pt\hbox{\pixel\vrule height0.125pt width12pt depth0.125pt\pixel}\kern-13.5pt
\raise3.25pt\hbox{\pixel\vrule height0.125pt width12pt depth0.125pt\pixel}\,}
\def\bs{\backslash}
\def\Colon{\colon\;\;}

{\bigbold Computation in Coxeter groups II. Minimal roots}

\bigskip

Bill Casselman\hfill\break\null
Mathematics Department\hfill\break\null
University of British Columbia\hfill\break\null
Canada\hfill\break\null
{\tt cass\char64math.ubc.ca}

\bigskip

{\bold Abstract}.  {\sl In the recent paper (Casselman, 2001)
I described how a number of ideas
due to Fokko du Cloux and myself could be incorporated into
a reasonably efficient program to carry out multiplication in
arbitrary Coxeter groups.  At the end of that paper
I discussed how this algorithm could be used to
build the reflection table of minimal roots,
which could in turn form the basis of
a much more efficient multiplication algorithm.
In this paper, following a suggestion of Bob Howlett,
I explain how results due to Brigitte Brink
can be used to construct the minimal root reflection table
directly and more efficiently.}

\intro
A {\bold Coxeter system} is a pair $(W, S)$ where
$W$ is a group with a set of generators $S$ 
and relations
$$ (st)^{m_{s,t}} = 1 $$
for pairs $s$, $t$ in $S$.
These groups play an important role in mathematics far from visible
in this simple definition.  Much work has been done on them,
particularly on finite and affine Weyl groups, but
many phenomena involving Coxeter groups
remain unexplained, and it is likely that computer
explorations will be even more significant in the future than
they have been so far.
These groups become extremely complex as the size of $S$ grows, however, 
and computer programs dealing with them must be extremely efficient
to be useful.  Conventional symbolic algebra packages normally fail
to handle the difficulties satisfactorily.   This paper is the second of a
series in which I describe what may be even in the long term
the most efficient algorithms to do basic computations in arbitrary Coxeter groups,
at least on serial machines.  These algorithms depend strongly on
mathematical results of Brigitte Brink, Robert Howlett, and Fokko du Cloux.

The first problem one encounters is how to multiply two elements of the group,
which is by no means a simple matter.
For familiar examples such as $S_{n}$ 
there are many satisfactory solutions, and for the Weyl groups of
Kac-Moody algebras one can use the representation of the group
on the root lattice to reduce many problems to integer arithmetic.
But arbitrary groups, even sometimes
the extraordinary finite groups $H_{3}$ and $H_{4}$,
are computationally more demanding.  Even with Weyl groups
integer overflow may occur when working with the root lattice.
In addition, solving computational
problems in this domain leads to mathematically interesting questions as well.

The best combinatorial solutions---certainly, those
of greatest theoretical interest---seem to be those in which
elements of the group are represented in terms of products of elements of $S$.
More explicitly, impose an ordering on $S$.
Suppose $w$ to be equal to the product
$s_{1}s_{2} \dots s_{n}$ as a product of generators in $S$,
and for each $i$ let $w_{i}$ be the partial product  $s_{1} \dots s_{i}$.  
This expression for $w$ is called its
\ISL/
{\bold normal form} if for each $i$ the element $s_{i}$ is least in $S$ such 
that $\ell(w_{i} s_{i}) < \ell(w_{i})$.   It is
the {\bold short}est representation of $w$ least in {\bold lex}icographical
order when read {\bold in} re{\bold verse}.   The multiplication
problem can now be posed: {\sl given the normal form
of $w$, what is that of $sw$?}  The first general solution was implicit
in work published by Jacques Tits in the mid 1960's.  
It had the virtue of
being entirely combinatorial in nature, working only
with the normal form of $w$, but it was
extremely inefficient, requiring both time and storage space roughly
exponentially proportional, or worse, to the length $n$ of the normal form of $w$.

A second solution was proposed by Fokko du Cloux around 1990, but explicitly 
described only for finite Coxeter groups.  It was this 
that I worked out in more detail in (Casselman, 2001).  It is
mildly recursive in nature, of unknown theoretical complexity but in practice
not unreasonably slow.  It also has the tremendous virtue of requiring little 
machine memory.
If one wants to be able to carry out just a few multiplications, it is likely
that this method cannot be improved on.  But if one wants to
carry out a great many multiplications for a given Coxeter group,
then it is probably best to build an auxiliary data structure that will
make multiplication much simpler.  A structure perfectly suited to
this purpose is a table
describing how the generators of $S$ act on
the {\bold minimal roots} of (Brink and Howlett, 1993).
As will be recalled later on,
it can be used in a multiplication algorithm of complexity essentially
linear in the length of $w$, allowing multiplication to take
place without backtracking in a left-to-right scan of the normal
form of $w$.

In the paper (Casselman, 2001) I proposed that du Cloux's
algorithm could be used as a kind of bootstrap
to build this reflection table, but as I mentioned 
at the end of that paper it is not difficult to find groups
where the number of minimal roots is rather small,
hence the table itself also fairly small, but the task
of building it in this way extremely tedious.
This is unsatisfactory.

After I explained this method to Robert Howlett, he and I experimented 
with other tools, and after a short while he
suggested that the results of Brigitte Brink's thesis
seemed to be applicable.  In this thesis, almost all of which
has been summarized in her published articles, she 
describes rather explicitly
what the set of minimal roots looks like for any Coxeter group.
However, it was not at all apparent at first glance, at least to me,
that what she does is computationally practical.
Howlett showed that my first impression was false by 
producing a program in the symbolic algebra package {\sl Magma}
that builds the minimal root
reflection table of an arbitrary Coxeter group.
This implementation was much slower than subsequent programs
written by Fokko du Cloux and myself, but already dealt easily
with a previously difficult group described in my paper.
It might be worthwhile to say that
the new algorithm depends on several crucial results
of Brink's thesis, but not on her more or less
explicit lists, interesting though they may be.

The principal drawback of
the new algorithm as the basis of multiplication
is that it claims an amount of machine memory proportional
to the number of minimal roots.  
The number of minimal roots can grow quite fast with the size of $S$,
so this is not a negligible objection.   It thus exemplifies the usual programming 
trade-off between time and space.
Familar Coxeter groups are
deceptive in this regard---for classical series of finite and affine groups
the number of minimal
roots is roughly $|S|^{2}$, whereas one of the first really
difficult groups Howlett's program dealt with
was an exotic Coxeter group with $|S| = 22$ 
and several hundred thousand minimal roots.
(Howlett's {\sl Magma} program took 69 minutes to handle it, he tells me,
whereas a program in {\sl C} written by Fokko du Cloux took 10 seconds 
on a comparable machine.)
Of course, such groups are so complicated that there is not
much exploration one can expect to do with them anyway, and in practice the
number of minimal roots is only a relatively minor impediment.

In the rest of this paper I will first recall the basic properties of minimal roots;
explain how they are to be used in multiplication; summarize
results of Brigitte Brink's thesis in the manner most
useful for my purposes; and finally explain how
to apply these results in practice to construct the minimal root reflection table.
I should say that what I am going to come up with is somewhat different 
from Howlett's program, and also quite
different from a related program of Fokko du Cloux.
It uses a few more advanced
features of Brink's thesis than Howlett's program.
It is less sophisticated than du Cloux's program,
but has the advantage of being possibly clearer in conception.   

There is not much originality in this paper.
Once Howlett had suggested using the results of Brink's 
thesis, the path to a practical program was probably almost determined.
On the other hand, this path did not appear obvious to me
until after a great deal of experiment, and I feel it will be useful
to place here a record of the outcome.   I should also
say that the order of exposition and its occasional geometric emphasis are new. 
I will be satisfied if this exposition awakens interest in the work of
Brink, Howlett, and du Cloux.
This whole area of research is one of great charm, but
it has apparently and unfortunately awakened little interest at large.
This is likely due to the fact that it lies somewhere
in the great no man's land between 
pure mathematics and practical computation---{\sl terra incognita}
on most maps, even in the 21st century.
I must also record that both my programs and this exposition have
benefited enormously from conversations with Fokko du Cloux.
As I learned long ago to my own chagrin,
competing directly in programming with 
him can be an embarrassing
experience, but working alongside him is, on the contrary, at once
pleasant and educational.


\Section{Minimal roots}
The algorithms to be described are intimately related to 
geometric realizations of Coxeter groups.
Following Brink and Howlett, I shall
work here exclusively with the {\bold standard realization}.
To each element $s$ of $S$ is associated an element
$\alpha_{s}$ of a basis $\Delta$ of a real vector space $V$.
An inner product is defined on $V$ according to the formula
$$ \alpha_{s} \Dot \alpha_{t} = -\cos (\pi/m_{s,t}) 
	= - \big(\zeta^{\phantom{-1}}_{2m} + \zeta^{-1}_{2m}\big)/2 \quad \big(m = m_{s, t},
\;  \zeta_{2m} = e^{2 \pi i / 2m}\big) 
\; . $$
In particular, $\alpha_{s} \Dot \alpha_{s} = 1$,
$\alpha_{s} \Dot \alpha_{t} = -1$ whenever
$s$ and $t$ generate an infinite group, and
$-1 < \alpha_{s} \Dot \alpha_{t} \le 1/2$ whenever they generate a finite
non-abelian group.
A representation of $W$ on $V$
is determined by mapping elements of $S$ to reflections:
$$ s v = v - 2 \, ( v \Dot \alpha_{s} ) \alpha_{s} \; . $$
Recall that the {\bold Coxeter graph}
has as nodes the elements of $S$, and an edge linking $s$ and $t$ if $m_{s,t}$,
which is called the {\bold degree} of the edge, is more than $2$.
In other words, the edges link $s$ and $t$ with $\alpha_{s} \Dot \alpha_{t} \ne 0$.
In case $\alpha \Dot \beta \ne 0$, I write $\alpha \sim \beta$.  
A {\bold simple link} in the graph is one with degree $3$, a
{\bold multiple link} one of higher degree.
The formula
for reflection amounts to the specification that if\/ 
$\lambda = \sum \lambda_{\gamma} \gamma$ then under reflection by $s_{\alpha}$ 
all the coordinates of $\lambda$ remain constant 
except that indexed by $\alpha$, which changes to
$$ (s\lambda)_{\alpha} = -\lambda_{\alpha} + \sum_{\alpha \sim \beta} [ -2 \, \alpha \Dot \beta ] \, \lambda_{\beta} \; . $$
In calculations, it is frequently best to work with $2 \, \alpha \Dot \beta$
rather than the dot-product itself,
since it is at once simpler and more efficient to do so.

This representation is faithful.
Let 
$$ C = \{ v \mid v \Dot \alpha_{s} > 0 \} \; .$$
The faces of this simplicial cone are the closures
of the open simplices
$$ C_{T} = \{ v \in V \mid \< \alpha, v > = 0 \hbox{ for } \alpha \in T,
	\<\alpha, v> > 0 \hbox{ for } \alpha \notin T \} $$
where $T \subseteq S$.  Thus $C$ itself is $C_{\emptyset}$,
and a face of codimension one is some $C_{\alpha}$.
The union of the domains
$w\overline{C}$ is a convex conical region called the {\bold Tits cone} ${\cal C}$
of the realization.  
The group $W$ acts discretely on ${\cal C}$, and the closure 
$\overline{C}$ of the open simplicial cone
is a fundamental domain for this action.  Each face of $w\overline{C}$ 
is the transform of a unique $\overline{C}_{T}$;
it is said to be {\bold labeled} by $T$.

A transform $w\alpha_{s}$ 
of an element of $\Delta$ is called
a {\bold root} of this realization.  
The root hyperplane $\lambda = 0$ is fixed by the
{\bold root reflection} $s_{\lambda}$ where $s_{w\alpha} = ws_{\alpha}w^{-1}$,
the conjugate of an element in $S$.
For any root $\lambda$, define 
$$ \eqalign {
{\cal C}_{\lambda = 0} &= \{ v \in {\cal C} \mid \<\lambda, v> = 0 \} \cr
{\cal C}_{\lambda \ge 0} &= \{ v \in {\cal C} \mid \<\lambda, v> \ge 0 \} \cr
{\cal C}_{\lambda \le 0} &= \{ v \in {\cal C} \mid \<\lambda, v> \le 0 \} \;. \cr
} $$

A root $\lambda$ is called {\bold positive} if
$\lambda > 0$ on $C$, {\bold negative} if $\lambda < 0$ on $C$.
Every root is either positive or negative,
which means that no root hyperplane  
ever intersects $C$.    A root is positive if and only if it
can be expressed $\lambda = \sum_{\alpha \in \Delta} \lambda_{\alpha} \alpha$
with $\lambda_{\alpha} \ge 0$.  

\satz{Proposition}{J}{Suppose that $\lambda$ and $\mu$ are positive roots.
The following are equivalent:
\smallbreak
\beginitems{3em}
\item{(a)}
$\lambda \Dot \mu = \cos(\pi k/m_{s,t})$ with
$0 < k < m_{s,t}$ for some $s$, $t$ in $S$;
\item{(b)} 
$|\lambda \Dot \mu| < 1$;
\item{(c)}
the group generated by the root reflections $s_{\lambda}$ and
$s_{\mu}$ is finite;
\item{(d)}
the hyperplanes 
$\lambda = 0$
and 
$\mu = 0$ intersect in the interior of ${\cal C}$.
\enditems
}
\endsatz

\proof/.  This is implicit in (Vinberg, 1971).
That (a) implies (b) is trivial; that (b) implies (c) is
a consequence of the simplest part of his discussion
of pairs of reflections.
If $G$ is any finite subgroup of $W$ then
the $G$-orbit of any point in the interior
of ${\cal C}$ will
be finite, and since ${\cal C}$ is convex the centre of mass
of the orbit will be in the interior of ${\cal C}$ as well.  This shows that (c)
implies (d).   Suppose $L = \{ \lambda = 0\} \cap \{ \mu = 0\} $
to contain a point in the interior of ${\cal C}$.
The whole of ${\cal C}$ is tiled by simplices
$wC_{T}$.  Points in the open cones $wC$ 
are fixed only by $1$ in $W$,
and those in some $wC_{\alpha}$ are fixed only
by $1$ and a single reflection.  Therefore
the linear space $L$ must contain some $wC_{T}$
with $|T| = 2$.  This means that $s$ and $t$ can be simultaneously
conjugated into $W_{T}$, and shows that (d) implies (a).
\qed

\satz{Corollary}{W}{If $\lambda$ and $\mu$ are distinct positive roots 
such that $|\lambda \Dot \mu| \ge 1$, then either
${\cal C}_{\lambda \ge 0}$ contains ${\cal C}_{\mu \ge 0}$
or 
${\cal C}_{\mu \ge 0}$ contains ${\cal C}_{\lambda \ge 0}$.}
\endsatz

\proof/.  If the intersection
of $\lambda = 0$ and $\mu = 0$ with the interior of
${\cal C}$ is null, then the sign of $\mu$ on the intersection of
$\lambda = 0$ with that interior will be constant.
Suppose it is negative.  I claim that the region
${\cal C}_{\mu \ge 0}$ is contained in the region
${\cal C}_{\lambda \ge 0}$.  If not, suppose that 
$P$ is a point of ${\cal C}$ with $\<\mu, P> \ge 0$
and $\<\lambda, P> < 0$.  If $P_{0}$ is
any point in the interior of $C$, then the open line segment $(P_{0}, P)$
is contained in the interior of ${\cal C}$, it is also contained in
the region ${\cal C}_{\mu \ge 0}$, and it
must intersect the hyperplane $\lambda = 0$.  At that point $\mu$ will be non-negative,
a contradiction.  If the sign is positive, switch $\lambda$ and $\mu$.\qed

Following Brink
and Howlett, I say that
the positive root $\lambda$ 
{\bold dominates} the positive root $\mu$ if the
region ${\cal C}_{\lambda \ge 0}$ 
 contains the region ${\cal C}_{\mu \ge 0}$. 
Another way of putting this is to require that 
$\lambda > 0$ in every
chamber $wC$ on which $\mu > 0$.
Loosely speaking, a {\bold minimal root} is one
which dominates only itself.  It is an important and remarkable
fact about arbitrary Coxeter
groups, first proved by Brink and Howlett, that the number of minimal
roots is always finite.  The significance of this
result does not appear in familiar cases,
since all roots are minimal
when $W$ is finite, and the minimal
roots for affine Coxeter groups
are rather simple in nature.
\medbreak

\leavevmode
\null
\hfill
\epsfxsize=1.6truein
\epsfbox{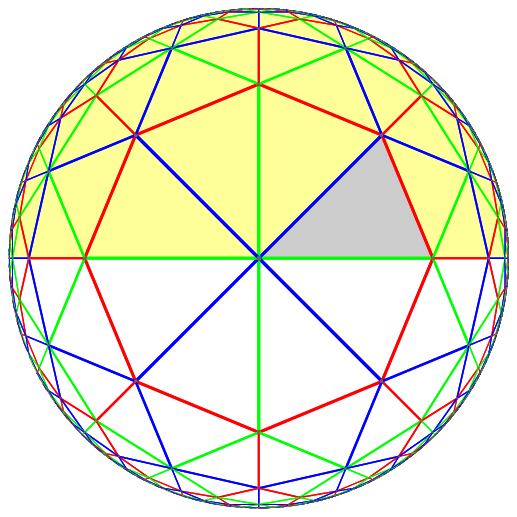}
\epsfxsize=1.6truein
\epsfbox{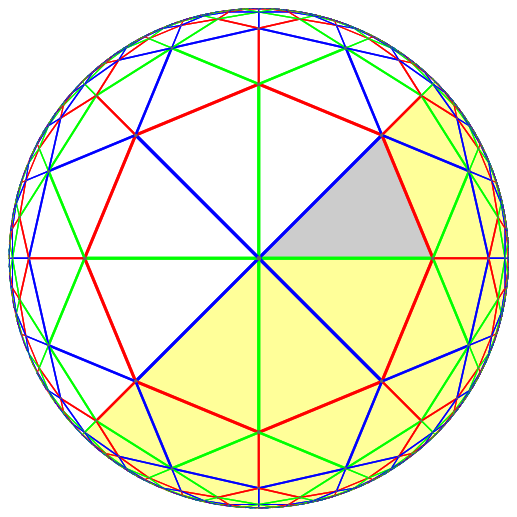}
\epsfxsize=1.6truein
\epsfbox{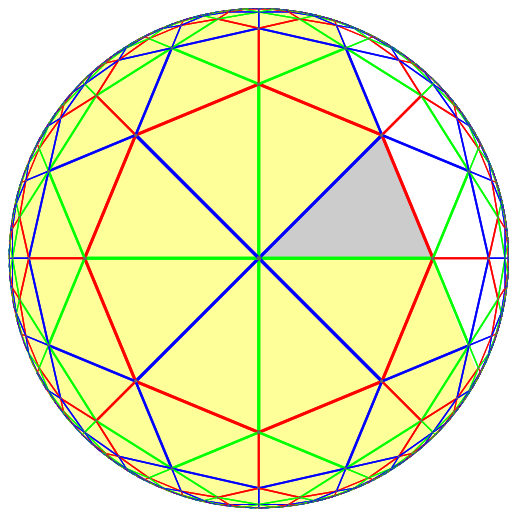}
\hfill
\null

\leavevmode
\null
\hfill
\epsfxsize=1.6truein
\epsfbox{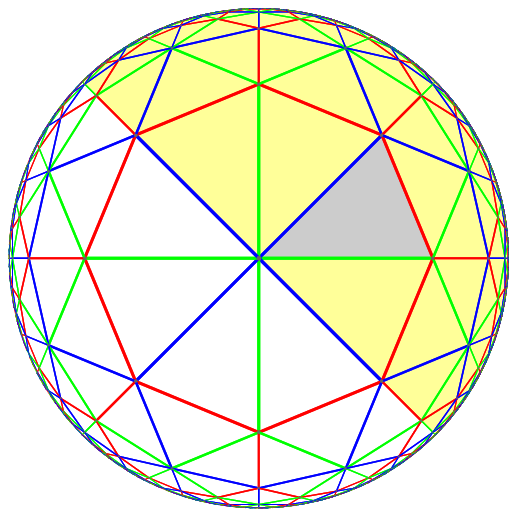}
\epsfxsize=1.6truein
\epsfbox{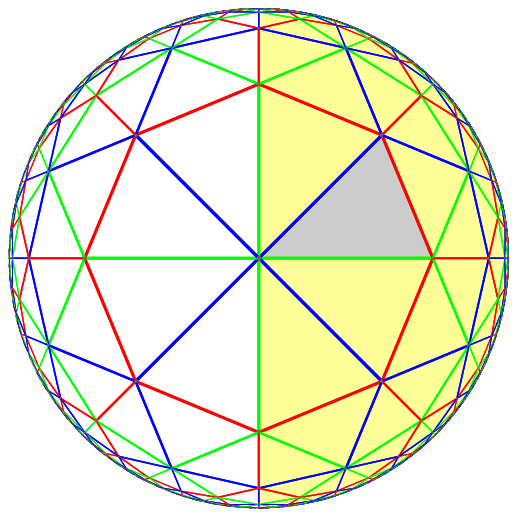}
\epsfxsize=1.6truein
\epsfbox{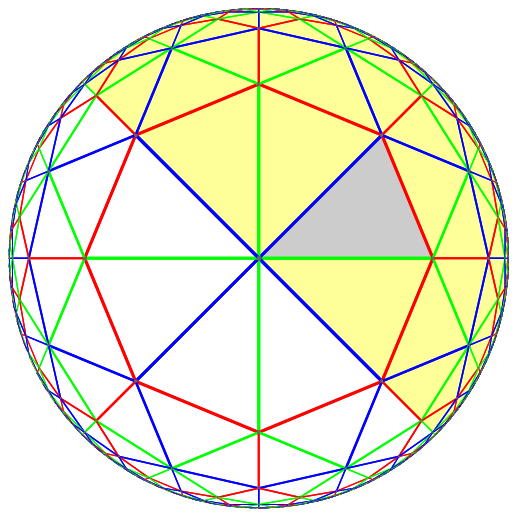}
\hfill
\null

\leavevmode
\null
\hfill
\epsfxsize=1.6truein
\epsfbox{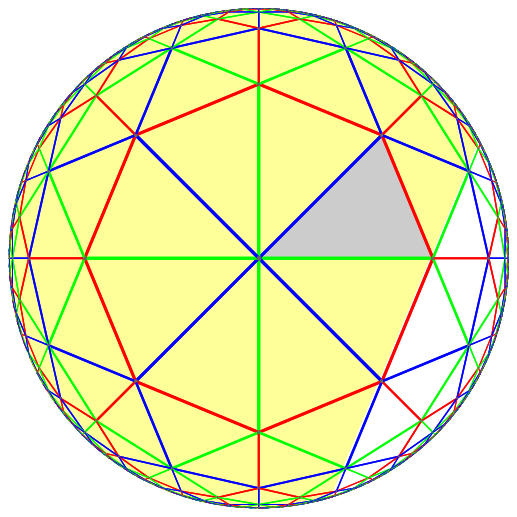}
\hskip0.4truein
\raise16pt
\vbox{\hsize=2.6truein
\overfullrule=0pt 
{\bold Figure 1.}
The domains $\lambda \ge 0$ for the seven minimal
roots of a Coxeter group with matrix
$$ \left[
\matrix{
1 & 3 & 4 \cr
3 & 1 & 3 \cr
4 & 3 & 1 \cr
}
\right]$$
(Exhibited in a slice through the Tits cone.)
}\hskip0.4truein
\hfill
\null

In the rest of this section I recall from
(Brink and Howlett, 1993) some basic facts about minimal roots.  
Proofs and statements will often be the same as theirs, 
and I'll refer to them when this happens. But in some places what I take
from them is only implicit in their paper, and in other places I
have given both statements and proofs a geometrical
flavour, somewhat as in (Casselman, 1995).

\satz{Lemma}{A}{{\rm (B \& H 2.2(ii))} Suppose $\lambda$ and $\mu$ to be positive
roots with $\lambda$ dominating $\mu$.  If $w\mu$ is positive, then 
$w\lambda$ is also positive, and dominates $w\mu$.}
\endsatz

\proof/.  Under the hypothesis, ${\cal C}_{w\lambda \ge 0}$ contains
${\cal C}_{w\mu \ge 0}$.  If $w\mu$ is positive, the second contains $C$, and therefore
$w\lambda > 0$ on $C$ as well.\qed

It is well known that $s$ in $S$ permutes the complement of $\alpha_{s}$ 
in the set of all positive roots, or in other words that if $s$ is in $S$ and $\lambda$
a positive root, then either $\lambda = \alpha_{s}$ and $s\lambda < 0$, or $s\lambda > 0$.
For the minimal roots, there is a similar range of options:

\satz{Proposition}{K}{If $\lambda$ is a minimal root and $s$ in $S$ then exactly one of these
three options holds:
\smallbreak
\beginitems{3em}
\item{$\bullet$} $\lambda = \alpha_{s}$ and $s\lambda < 0$;
\item{$\bullet$} $s\lambda$ is again minimal;
\item{$\bullet$} $s\lambda$ dominates $\alpha_{s}$.
\enditems}
\endsatz

\proof/.  If $\lambda \ne \alpha > 0$, then $s_{\alpha}\lambda > 0$.
If it is not minimal, say it dominates $\beta > 0$ with
$s_{\alpha}\lambda \ne \beta$.
If $\beta \ne \alpha$, then $s_{\alpha}\beta > 0$, and by the previous Lemma
$\lambda$ dominates $s_{\alpha}\beta$.  Since $\lambda$ is minimal,
this can only happen if $\lambda = s_{\alpha}\beta$, a contradiction.\qed

The {\bold minimal root reflection table} is a matrix $\rho$
of size $|S| \times N$, where $N$ is the number of minimal
roots.  The entries are either minimal
roots, or {\bold virtual} minimal roots
I arbitrarily label as $\ominus$ and $\oplus$.
I set $\rho(s, \lambda) = \mu$
if $s\lambda = \mu$ is minimal,
$\rho(s, \lambda) = \ominus$ if $\lambda = \alpha_{s}$,
and $\rho(s, \lambda) = \oplus$ if $s\lambda$ is a non-minimal
positive root.
This table is one of the fundamental data structures of
a Coxeter group.

The {\bold depth} $\delta(\lambda)$ of a positive root $\lambda$ is the length
of the smallest $w$ such that $w^{-1}\lambda < 0$.  
Equivalently, it is one less than the length of the shortest gallery
$C = C_{0}$, $C_{1}$, \dots , $C_{n}$ with $C_{n}$ in the region
$\lambda < 0$.  Recall that a gallery is a chain of chambers
with successive chambers sharing a common wall.  The depth is also
the smallest $n$ with 
$$ \lambda = s_{1}s_{2} \dots s_{n-1}\alpha $$
for some $\alpha \in \Delta$.  The elements of $\Delta$ itself,
for example, have depth $1$. 
A partial order is thus induced on the set
of all roots: $\lambda \preceq \mu$ if $\mu = w \lambda$ where
$\delta(\mu) = \ell(w) + \delta(\lambda)$.  This order gives rise to
the {\bold root graph}, whose nodes are positive roots,
 with edges $\lambda \rightarrow s\lambda$ whenever
$\lambda \prec s\lambda$.

\satz{Lemma}{B}{If $\lambda > 0$ and $s_{\alpha} \lambda$ dominates $\alpha$ 
then $\lambda \prec s_{\alpha}\lambda$.}
\endsatz

\proof/.  If $(C_{i})$ is a gallery from $C$ to the region where $s_{\alpha}\lambda < 0$,
there will be some pair $C_{i}$, $C_{i+1}$ sharing a wall on the hyperplane
$\alpha = 0$.  The reflected gallery will be shorter.  
Hence $\delta(\lambda) < \delta(s_{\alpha}\lambda)$.\qed

The minimal roots lie at the bottom of the root graph, in the sense that:

\satz{Proposition}{L}{{\rm (B \& H 2.2(iii))} If $\lambda \preceq \mu$ and $\mu$
is minimal, so is $\lambda$.}
\endsatz

\proof/.  This reduces to the case $\mu = s_{\alpha}\lambda$.
If $\mu$ is minimal but $\mu \ne \alpha$,
then either $s_{\alpha}\mu = \lambda$ is minimal,
or it dominates $\alpha$.  But by the preceding Lemma,
in the second case $\delta(\mu) < \delta(s_{\alpha}\mu)$,
a contradiction.\qed

\satz{Lemma}{C}{{\rm (B \& H 2.3) } If $\lambda > 0$ and $\lambda \Dot \alpha > 0$
then $s_{\alpha}\lambda \prec \lambda$.}
\endsatz

\proof/.  Suppose $w^{-1}\lambda < 0$
with $\ell(w) =\delta(\lambda)$.
If $w^{-1}\alpha < 0$ then we can write 
$$ \eqalign {
w &= s_{\alpha}s_{2} \dots s_{n} \cr
\lambda &= s_{\alpha}s_{2} \dots s_{n-1}\alpha_{n} \cr
s_{\alpha}\lambda &= s_{2} \dots s_{n-1}\alpha_{n} \cr
\delta(s_{\alpha}\lambda) = n-1 &< \delta(\lambda) = n \; . \cr
} $$

So we may assume that $w^{-1}\lambda < 0$ but
$w^{-1}\alpha > 0$.  Let $w$ have the reduced expression
$w = s_{1} \dots s_{n}$, and let $u = w s_{n} = s_{1} \dots s_{n-1}$.
Thus $w\alpha_{n} < 0$.
I claim that $u^{-1}s_{\alpha}\lambda < 0$, which implies that
$\delta(s_{\alpha}\lambda) = n-1$.

We can calculate that
$$ \eqalign {
	w^{-1}s_{\alpha}\lambda 
		&= w^{-1}\lambda - 2\, (\lambda \Dot \alpha) w^{-1}\alpha \cr
	u^{-1}s_{\alpha}\lambda 
		&= u^{-1}\lambda - 2\, (\lambda \Dot \alpha) s_{n} w^{-1}\alpha \; . \cr
} $$
From the first equation, we see that $w^{-1}\lambda$ 
will be negative under the assumption that $\lambda \Dot \alpha > 0$.
Since the depth of $\lambda$ is $n$, $\mu = u^{-1}\lambda$ will be positive.
But then $\mu > 0$, $s_{n} \mu < 0$ implies that $\mu = \alpha_{n}$.
So we have
$$	u^{-1}s_{\alpha}\lambda 
		= \alpha_{n} - 2\, (\lambda \Dot \alpha) s_{n} w^{-1}\alpha \; . $$
The root
$$ \alpha_{n} - 2\, (\lambda \Dot \alpha) s_{n} w^{-1}\alpha $$
is either positive or negative. 
In the second term, the positive root
$w^{-1}\alpha$ cannot be $\alpha_{n}$, since then $w\alpha_{n} < 0$.
Therefore
$s_{n}w^{-1}\alpha$ is a positive root also not $\alpha_{n}$,
and the second term is a negative root not equal to a multiple of $-\alpha_{n}$.
The sum has to be negative.\qed

\satz{Corollary}{X}{{\rm (B \& H 1.7)} Whenever $\lambda$ is a positive root and $s$ in $S$ there are
three possibilities: (a) $s\lambda \rightarrow \lambda$, (b) $s\lambda = \lambda$,
or (c) $\lambda \rightarrow s\lambda$, depending on
whether $\lambda \Dot \alpha_{s} > 0$, $\lambda \Dot \alpha_{s} = 0$,
or $\lambda \Dot \alpha_{s} < 0$, respectively.}
\endsatz

\satz{Corollary}{Y}{{\rm (B \& H 1.8)}  If\/ 
$\lambda = \sum \lambda_{\alpha} \alpha \preceq \mu = \sum \mu_{\alpha} \alpha$
then $\lambda_{\alpha} \le \mu_{\alpha}$.}
\endsatz

These results will soon be used to outline an algorithm to construct
the minimal root reflection table.
Processing roots will take place in order of increasing depth, applying this criterion:
if $\lambda \ne \alpha$ is a minimal root,
then $s_{\alpha}\lambda$ is no longer a minimal root if and only 
if $\lambda \Dot \alpha_{s} \le -1$.  Details follow, after I explain
in the next section how to use the minimal
root reflection table in multiplication.


\Section{Minimal roots and multiplication}
Suppose that $w$ has the normal form
$w = s_{1} \dots s_{n}$.  Then the normal
form of $sw$ will be obtained from that 
of $w$ by {\bold insertion} or {\bold deletion}:
$$ NF(sw) = s_{1} \dots s_{i} t s_{i+1} \dots s_{n} $$
or
$$ NF(sw) = s_{1} \dots s_{i-1} s_{i+1} \dots s_{n} \; . $$
How do we find where the insertion or deletion occurs,
and which of the two it is?  If insertion, what is inserted?

The algorithm explained in (Casselman, 2001)
handles these problems in geometric terms.
I define the {\bold \ISL/ tree} to be the directed graph
whose nodes are elements of $W$, with 
a link from $x$ to $y$ if $NF(y) = NF(x) \cat/ s$.
Its root is the identity element.
Each $y \ne 1$ in $W$ is the target of exactly one link,
labeled by the $s$ in $S$
least with $ys < y$.
Paths in this tree starting at $1$ 
match the normal forms of elements of $W$.
\medbreak

\captionedFigure{3truein}{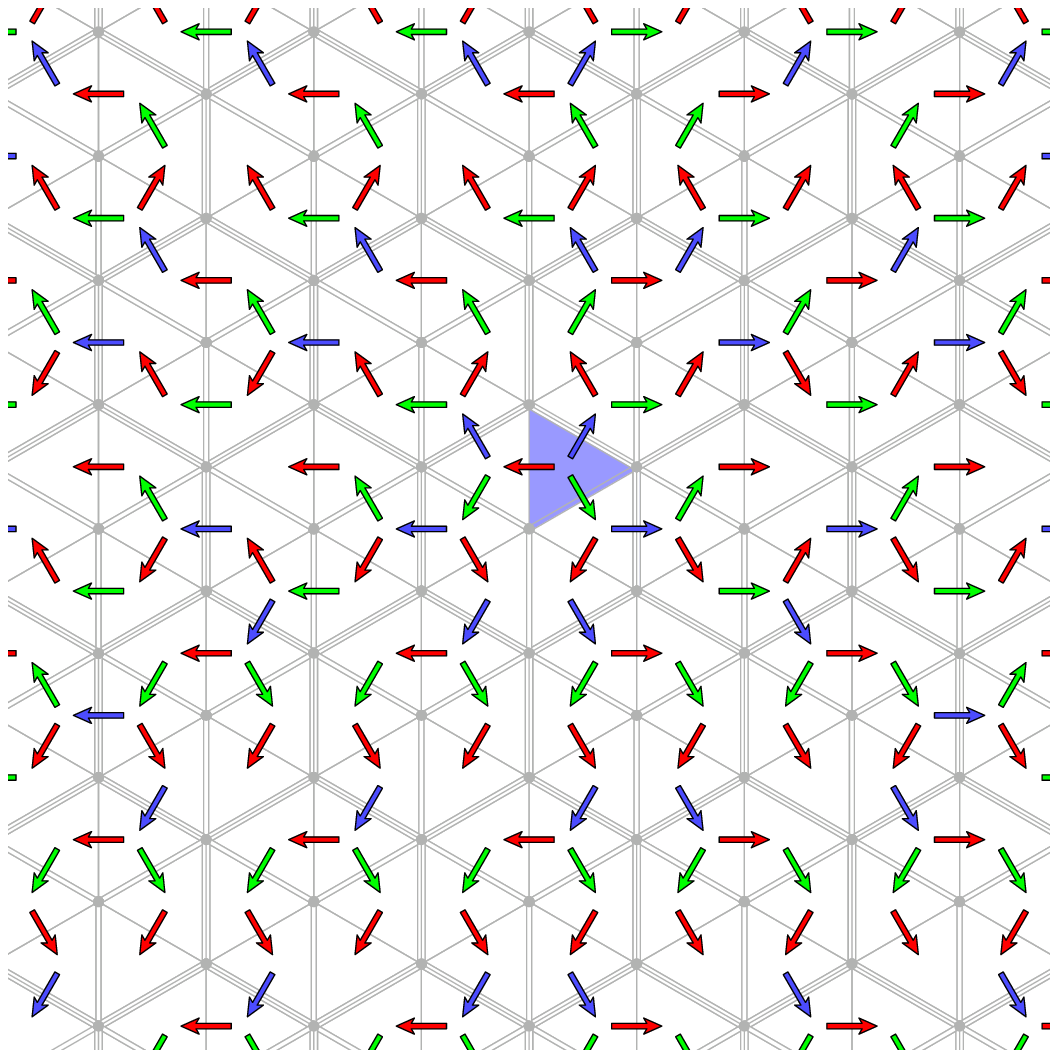}{The \ISL/ tree for $\widetilde{A}_{2}$.}

This tree gives rise to a geometrical figure in which
we put a link from $xC$ to $yC$ if there is an edge in the tree
from $x$ to $y$. 
The multiplication problem can now be formulated:
{\sl given an \ISL/ path from $C$ to $wC$, how can we find the \ISL/ path
from $C$ to $swC$?}  As explained in (Casselman, 2001),
the \ISL/ figure is very close to being symmetric
with reflection to the reflection $s$; only
links near the hyperplane $\alpha_{s} = 0$ are not
preserved upon reflection.  The chambers principally
affected are what I call {\bold exchange sites} of the \ISL/ tree---$wC$
is an exchange site if the chamber $wC$ has a face
on the hyperplane $\alpha_{s} = 0$, labeled by $t$ in $S$
which is less than the labelling of the \ISL/ link
entering $wC$.
\medbreak

\captionedFigure{3truein}{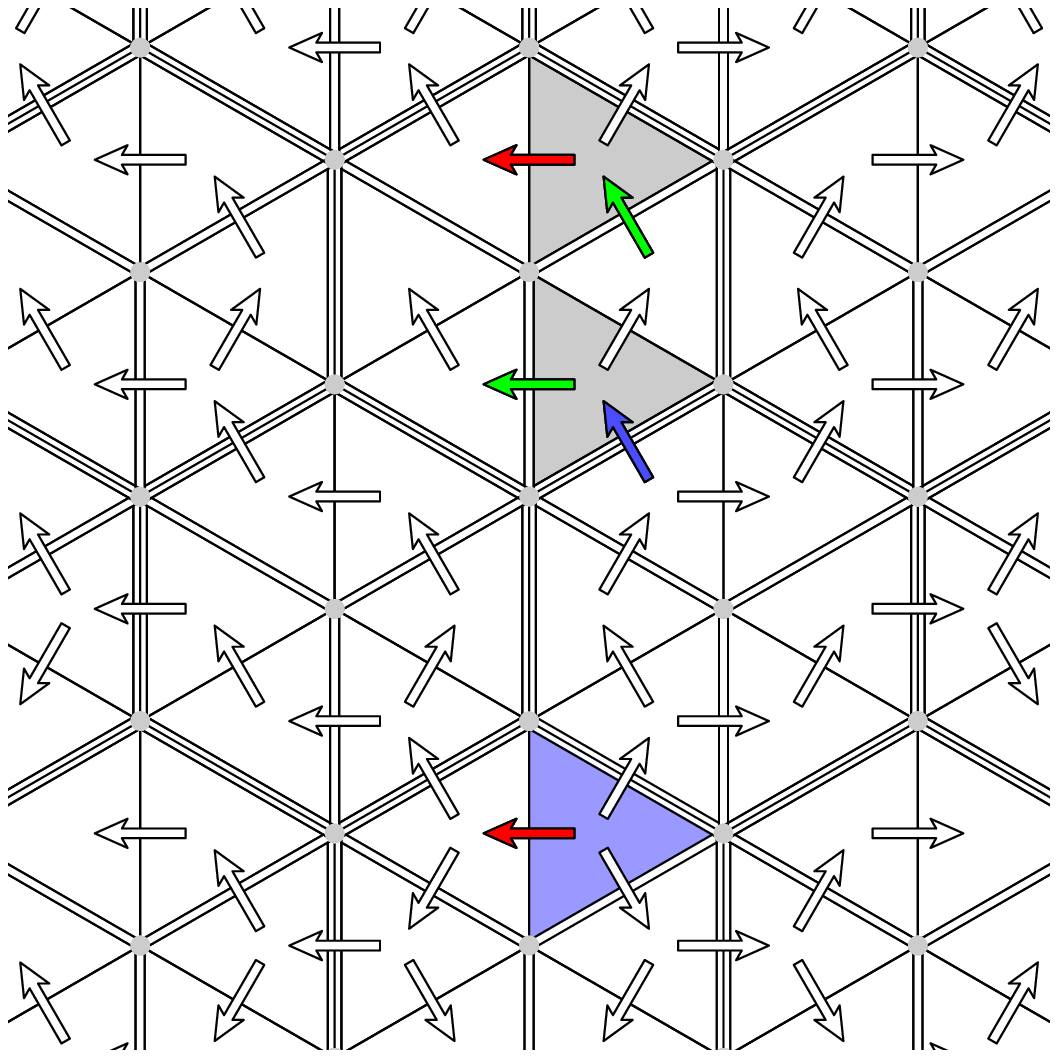}%
{Some of the 
exchange nodes for reflection by $\[1]$ in $\widetilde{A}_{2}$.  Links
modified under reflection are colored.  
The fundamental chamber is at the bottom.}

The relationship between this notion and
multiplication is simple.  Suppose $NF(w) = s_{1} \dots s_{n}$ and for each $i$
let $w_{i}$ be the partial product $s_{1} \dots s_{i}$.  Suppose that $w_{k}C$
is the last exchange site in the gallery $C = C_{0}$, \dots, $C_{n} = w_{n}C$,
and suppose that $t$ labels the wall on $\alpha = 0$.  Then either the gallery
crosses $\alpha = 0$ there or it does not.  In the first case,
$t = s_{k+1}$ and 
$$ NF(sw) = s_{1} \dots s_{k} s_{k+2} \dots s_{n} \; . $$
In the second, 
$$ NF(sw) = s_{1} \dots s_{k} t s_{k+1} \dots s_{n} \; . $$

The question is now: {\sl How do we recognize \ISL/ exchange sites?}
This is where roots and, even more pleasantly, minimal
roots come in. 

\satz{Lemma}{D}{Suppose $w = s_{1} \dots s_{n}$.
The chamber $wC$ is an \ISL/
exchange site if and only if\/ 
$w^{-1} \alpha = \beta$ with $s_{\beta} < s_{n}$.} 
\endsatz

\proof/. The wall of $wC$ labeled by $\beta$
is embedded in the hyperplane $w\beta = 0$.\qed

This suggests the following algorithm to calculate $NF(sw)$:
Read the string $s_{1} \dots s_{n}$ from left to right.   
As we read, keep track of
an element $t$ in $S$, an index $k$, and a root $\lambda$.  
The index $k$ is that of the last exchange,
the element $t$ is the label of the last exchange, $\lambda = w^{-1}\lambda$.
To start with $k = -1$, $t = s$, $\lambda = \alpha_{s}$.
As we read $s_{i}$, $\lambda$ changes to $s_{i}\lambda$.
If $\lambda = \beta$ with $s_{\beta} < s_{i}$, $t$ changes to $s_{\beta}$
and $k$ to $i$.  When we are through scanning, we have the data determining
the last exchange.

So far, minimal roots have not come into play.
But in fact, in calculating the reflected root
$s\lambda$ we only need to apply the 
minimal root reflection table.  If $\lambda$ is minimal
and $s\lambda = \ominus$,
then the chambers $w_{i-1}C_{i-1}$ and $C_{i} = w_{i-1}s_{i}C_{i}$
lie on opposite sides of the hyperplane $\alpha = 0$,
so we know we are looking at a deletion.
If $s\lambda$ is minimal and $s\lambda = \oplus$ then
we know the gallery has crossed a hyperplane $\lambda = 0$ which separates
$wC$ from $\alpha = 0$.  The exchange we last recorded 
will be the last we ever record, because
the gallery cannot recross that hyperplane to reach $\alpha = 0$.

We now know how to use the minimal root reflection table to multiply,
and it remains to tell how to construct it.


\Section{A rough guide to constructing the minimal roots}
A straightforward procedure to list all minimal roots,
along with a description of how the elements
of $S$ act on them, is not difficult to sketch.

The data determining a minimal root 
$\lambda$ are an index,
a coordinate vector, and a list of reflections $s\lambda$.
There are at least two reasonable choices of coordinates.
The coordinates I use in this note are those in the expression 
of a root as a linear combination of the roots in $\Delta$:
$$ \lambda = \sum_{\alpha \in \Delta}\lambda_{\alpha} \alpha $$
The other reasonable choice is the array $(\lambda \Dot \alpha)$,
which du Cloux has discovered to be a choice somewhat more efficient 
if also more difficult to work with in programming.
Indices are assigned in the order in which the minimal
roots are first calculated; 
the roots in $\Delta$ are assigned the first $|S|$ integers.
Furthermore, in the preliminary version
of the algorithm we shall maintain
a look-up table of some kind, enabling us to find
a root's data given its coordinates.

A minimal root is to be dealt with in two stages:
(1) {\bold defining} and (2) {\bold finishing}.  In the definition,
a root is calculated, assigned an index, entered in the look-up table,
and put on a waiting list.  Items are taken off
this list in the order in which they are put
on---the list is, in programmers' terminology,
a FIFO list ({\bold f}irst {\bold i}n, {\bold f}irst {\bold o}ut) or {\bold queue}.

A root $\lambda$ is finished when it is taken off the queue
and all of the reflections $s\lambda$ not already known are found.
At this moment, all roots with depth less than that of
$\lambda$ have been finished, so that in particular we have recorded
all reflections $\mu = s\lambda \prec \lambda$.
We must calculate the $s\lambda$
not yet determined;
the problem is to figure out whether $s\lambda = \lambda$,
$s\lambda$ is a minimal root, or $s\lambda = \oplus$.
Of course we can tell by inspection
whether $s \lambda = \lambda$.   Otherwise we
know that $\lambda \prec \mu = s\lambda$.   
There are three possibilities:
(1) We have already 
defined $\mu$, in which case it will be registered in
our look-up table.  We just add to its record that $s\lambda = \mu$
and $s\mu = \lambda$.   
(2) The root $\mu$ is not minimal.
We can tell this when we calculate $\mu$, since in this case
$\lambda \Dot \alpha_{s} \le -1$.  We set $s\lambda = \oplus$.
(3) Otherwise $-1 < \lambda \Dot \alpha_{s} < 0$, and the root $\mu$
is a new minimal root not yet defined.
We do so by assigning it an index, putting it in the look-up table,
recording the reflection $s\lambda = \mu$, and putting $\mu$
in the process queue.   Then we go on to the next root in
the queue, as long as there is something there.

There are a number of problems 
that occur in trying to implement this is in practice.
The first is that maintaining the look-up table is
somewhat cumbersome.  As du Cloux first observed,
one of the results from Brink's
thesis allows us to dispense entirely it.
The following is a endsatz of the main Theorem of (Brink, 1995).

\satz{Lemma}{E}{Suppose $s \ne t$ in $S$.
Assume the group $W_{s,t}$ they generate
to be finite.  Suppose $\lambda$ to be any root whose 
$W_{s,t}$-orbit doesn't intersect $\Delta$.  
Then $W_{s,t}$ acts freely on this orbit.}
\endsatz

How can this be applied in our algorithm to find all minimal roots?
Suppose we are in the process of
defining a minimal root $\mu$, which we may
assume not to be dihedral.
We calculate it as $\mu = s\lambda$ with $\lambda \prec \mu$.  If $t$ is
an element of $S$ with $t\lambda \prec \lambda$,
then $W_{s,t}$ will necessarily be finite.  In the orbit of $\mu$ under
$W_{s,t}$ there will be a unique element $\nu$ of least depth,
equal to $w_{s,t} \mu$ where
$w_{s,t}$ is the unique longest element of the dihedral group
$W_{s,t}$.   We can calculate it by finding alternately $t\lambda$, $st\lambda$,
etc. until the chain starts to ascend in depth.   And the difference in depth
between $\mu$ and $\nu$ 
will be exactly $m_{s,t}$ precisely when $t\mu \prec \mu$.
The element $t\mu$ will then be
$(tw_{s,t}) \nu$, which can be calculated from the reflection tables
of elements of depth less than $\lambda$.  In other words, as soon
as $\mu$ is defined we can calculate all the {\bold descents}
$t\mu \prec \mu$---without doing any new calculations!   

In the algorithm above we can therefore eliminate reference to the look-up table.
The procedure sketched above requires, however,
that we know whether the reflection of
a root has less depth or not, so we add to
the data of a root $\lambda$ its {\bold descent set}, the set
of $s$ with $s\lambda \prec \lambda$.
When we define a root, before we put it in the queue, we calculate all of its descents.

Because the result above is
applicable only to
roots not in the orbit of $W_{s,t}$,
all the dihedral roots have to be handled somewhat specially---they
should all be defined immediately after the roots in $\Delta$,
and those of depth greater than $2$ have to be fed
into the queue at the right moment.

There is another simplification.
The following was observed
by Brink and Howlett, and used
as an important part of their proof
that the set of minimal roots is finite:

\satz{Proposition}{M}{If $s\lambda = \oplus$
and $\lambda \prec \mu$ then $s\mu = \oplus$ as well.}
\endsatz

\proof/.  For $\alpha$, $\beta$ in $\Delta$, $\lambda$ any root 
$$ s_{\alpha} \lambda \Dot \beta = \big(\lambda  - 2 (\lambda \Dot \alpha) \alpha \big) \Dot \beta
	=  (\lambda \Dot \beta) - 2 (\lambda \Dot \alpha)(\alpha \Dot \beta) $$
which means that for $\beta \ne \alpha$ the dot product
$\lambda \Dot \beta$ decreases under reflection by $\alpha$,
as long as $\lambda \prec s_{\alpha}\lambda$.
Therefore once we have calculated that $s\lambda = \oplus$ 
this remain true for all $\mu$ with $\lambda \prec \mu$.\qed

I call an $s$ with $s\lambda = \oplus$ a {\bold lock} at $s$.
This Proposition tells us that locks are inherited as we go
up the root graph and these inherited locks, like descents,  
can be immediately recorded upon definition.

Eliminating the look-up table is a great simplification, noting 
inherited locks
a minor one.  In simple cases---when none of
the Coxeter matrix entries $m_{s,t}$ are greater than $3$---there
is probably no better algorithm than the one
sketched above, but incorporating the simple modifications concerning descents
and locks.
In this case, all doubled dot products
$2\, (\lambda \Dot \alpha)$ will be integers. 
When $\lambda \prec s_{\alpha}\lambda$ this 
doubled dot-product will be a negative integer,
and if it is anything other than $-1$ the root $s\lambda$ will
no longer be minimal.  What could be simpler?

If there are entries $m_{s,t} > 3$ this still works,
at least in principle, although there are a few annoyances
that arise.
All dot products in the range
$(-1, 0)$ will be of the form $\cos \big(\pi k/m_{s,t}\big)$ for
$k < m_{s,t}/2$.  So in order to tell whether $s\lambda$ dominates
$\alpha_{s}$ we just have to compare the dot product
$\lambda \Dot \alpha_{s}$ to numbers in this finite set.  
There is one potentially nasty problem, however.
All we can assert {\sl a priori}
about the coordinates of
$\lambda$ is that they will be in the real cyclotomic field
generated by $\cos (2 \pi /M)$, where $M$ is the lowest common
denominator of the numbers $2 m_{s,t}$.
The amount of
work involved would seem to grow disturbingly as $M$ grows.
It is Brink's thesis that tells us that this apparent
difficulty is spurious.  In the next section I follow her
in examining more carefully the coefficients of
minimal roots.


\Section{Root coefficients}
We want to
know something about what possible values can occur for the 
coefficients of roots and, especially, minimal roots.
Roots are produced by a sequence of
reflections from the basic roots in $\Delta$.
As is the case in many algorithmic processes,
although we understand each single step---here, reflection---quite well,
the overall development is not so clear.

We begin with a discussion of dihedral roots.

Suppose $S$ to possess two elements $s$ and $t$,
corresponding to roots $\alpha$ and $\beta$.  Let $m = m_{s,t}$.
What are all the positive roots of the system?

Let $z = \zeta_{2m} = e^{\pi i/ m} = e^{2 \pi i/ 2 m}$.  Then
$2 \, \alpha \Dot \beta = 2 \cos(\pi /m) = (z + z^{-1})$ and
$$ \eqalign {
s \lambda &= \lambda - 2 (\lambda \Dot \alpha) \alpha \cr
 2 \, (p\alpha + q\beta) \Dot \alpha &= 2p - (z+z^{-1})q \cr
s (p\alpha + q\beta) &= \big[-p + (z+z^{-1})q\,\big] \alpha + q\beta \cr
t \lambda &= \lambda - 2 (\lambda \Dot \beta) \beta \cr
 2 \, (p\alpha + q\beta) \Dot \beta &= 2q - (z+z^{-1})p \cr
t (p\alpha + q\beta) &= p \alpha + \big[-q + (z+z^{-1}) p\, \big]\beta \; . \cr
} $$
Suppose $m = \infty$.  We get in succession (in tabular form)
as transforms of $\alpha$ 
\medbreak
\settabs\+ xxxxxxxxxx & xxxxxxxxxxxxxxx & xxxxxxxxxxxxxxx & xxxxxxxxxxxxxxx & \cr
\+ &  & \hfill coefficient of $\alpha$ & \hfill coefficient of $\beta$ & \cr
\+ & \hfill $\alpha$ & \hfill $1$ & \hfill $0$ & \cr
\+ & \hfill $t\alpha$ & \hfill $1$ & \hfill $2$ & \cr
\+ & \hfill $st\alpha$ & \hfill $3$ & \hfill $2$ & \cr
\+ & \hfill \dots & \hfill \dots & \hfill \dots & \cr
\+ & \hfill $(st)^{n}\alpha$ & \hfill $2n+1$ & \hfill $2n$ & \cr
\+ & \hfill $t(st)^{n}\alpha $ & \hfill $2n+1$ & \hfill $2n+2$ & \cr
\+ & \hfill \dots & \cr

and similarly for the transforms of $\beta$.  Of these, only
$\alpha$ and $\beta$ are minimal roots.

If $m=2$, the positive roots are $\alpha$ and $\beta$.

Now suppose $2 < m < \infty$.  Then we get

\medbreak
\settabs\+ xxxxxxxxxx & xxxxxxxxxxxxxxx & xxxxxxxxxxxxxxx & xxxxxxxxxxxxxxx & \cr
\+ &  & \hfill coefficient of $\alpha$ & \hfill coefficient of $\beta$ & \cr
\+ & \hfill $\alpha$ & \hfill $1$ & \hfill $0$ & \cr
\+ & \hfill $t\alpha$ & \hfill $1$ & \hfill $z+z^{-1}$ & \cr
\+ & \hfill $st\alpha$ & \hfill $z^2 + z^{-2} + 1$ & \hfill $z + z^{-1}$ & \cr
\+ & \hfill \dots & \hfill \dots & \hfill \dots & \cr

\+ & \hfill $(st)^{n}\alpha$ & \hfill $C_{n+1}$ & \hfill $C_{n}$ & \cr
\+ & \hfill $t(st)^{n}\alpha$ & \hfill $C_{n+1}$ & \hfill $C_{n+2}$ & \cr

where the $C_{n}$ are solutions of the difference equation
$C_{n+2} = -C_{n} + (z+z^{-1}) C_{n+1}$ with initial values
$C_{0} = 0$, $C_{1} = 1$, giving
$$ C_{n} = { z^{n} - z^{-n} \over z - z^{-1} }  
		= z^{n-1} + z^{n-3} + \cdots + z^{-(n-3)} + z^{-(n-1)} ={  \sin (\pi n/m) \over \sin(\pi/m) } \; . $$
Note that $C_{n} > C_{n-1}$ as long as the real part of $z^{n}$ is
positive, which happens as long as $n < m/2$.
There are $m$ roots in all for a dihedral system,
which agrees with this remark.

For $m=3$ we get the table
\medbreak
\settabs\+ xxxxxxxxxx & xxxxxxxxxxxxxxx & xxxxxxxxxxxxxxx & xxxxxxxxxxxxxxx & \cr
\+ &  & \hfill coefficient of $\alpha$ & \hfill coefficient of $\beta$ & \cr
\+ & \hfill $\alpha$ & \hfill $1$ & \hfill $0$ & \cr
\+ & \hfill $\beta$ & \hfill $0$ & \hfill $1$ & \cr
\+ & \hfill $t\alpha = s\beta$ & \hfill $1$ & \hfill $1$ & \cr

\medbreak
For $m=4$:
\medbreak
\settabs\+ xxxxxxxxxx & xxxxxxxxxxxxxxx & xxxxxxxxxxxxxxx & xxxxxxxxxxxxxxx & \cr
\+ &  & \hfill coefficient of $\alpha$ & \hfill coefficient of $\beta$ & \cr
\+ & \hfill $\alpha$ & \hfill $1$ & \hfill $0$ & \cr
\+ & \hfill $\beta$ & \hfill $0$ & \hfill $1$ & \cr
\+ & \hfill $t\alpha = st\alpha$ & \hfill $1$ & \hfill $\sqrt{2}$ & \cr
\+ & \hfill $s\beta = ts\beta$ & \hfill $\sqrt{2}$ & \hfill $1$ & \cr

\medbreak
For $m=5$, with $\theta = (1 + \sqrt{5})/2$:
\medbreak
\settabs\+ xxxxxxxxxx & xxxxxxxxxxxxxxx & xxxxxxxxxxxxxxx & xxxxxxxxxxxxxxx & \cr
\+ &  & \hfill coefficient of $\alpha$ & \hfill coefficient of $\beta$ & \cr
\+ & \hfill $\alpha$ & \hfill $1$ & \hfill $0$ & \cr
\+ & \hfill $\beta$ & \hfill $0$ & \hfill $1$ & \cr
\+ & \hfill $t\alpha$ & \hfill $\theta$ & \hfill $1$ & \cr
\+ & \hfill $s\beta$ & \hfill $1$ & \hfill $\theta$ & \cr
\+ & \hfill $st\alpha = ts\beta$ & \hfill $\theta$ & \hfill $\theta$ & \cr

\medbreak
For $m=6$:
\medbreak
\settabs\+ xxxxxxxxxx & xxxxxxxxxxxxxxx & xxxxxxxxxxxxxxx & xxxxxxxxxxxxxxx & \cr
\+ &  & \hfill coefficient of $\alpha$ & \hfill coefficient of $\beta$ & \cr
\+ & \hfill $\alpha$ & \hfill $1$ & \hfill $0$ & \cr
\+ & \hfill $\beta$ & \hfill $0$ & \hfill $1$ & \cr
\+ & \hfill $t\alpha = st\alpha$ & \hfill $1$ & \hfill $\sqrt{3}$ & \cr
\+ & \hfill $s\beta = ts\beta$ & \hfill $\sqrt{3}$ & \hfill $1$ & \cr
\+ & \hfill $st\alpha = tst\alpha$ & \hfill $2$ & \hfill $\sqrt{3}$ & \cr
\+ & \hfill $ts\beta = sts\beta$ & \hfill $\sqrt{3}$ & \hfill $2$ & \cr

\medbreak
The table for $m > 6$ looks somewhat similar to the last, except that it is longer.
All coefficients beyond the first $6$ rows are greater than $2$.

We can use these results, following Brink's thesis, to prove:

\satz{Proposition}{N}{The coefficients of
the positive roots are polynomials in the constants
$$ {\zeta_{2m}^{\,n\phantom{-}} - \zeta_{2m}^{-n} 
	\over \zeta^{\phantom{-1}}_{2m} - \zeta_{2m}^{-1} }\quad (n \le m/2) $$
with coefficients among
the natural numbers ${\Bbb N}$.}
\endsatz

\proof/.  By induction on the depth of a root.
If the depth is $1$ there is no problem.
Otherwise suppose 
$$ \lambda = w \alpha = s_{1} \dots s_{n-1}\alpha_{t} $$
with $\alpha$ in $\Delta$, $\delta(\lambda) = n$.
Let $s = s_{n}$.  Choose $y \ne 1$ maximal in
$W_{s,t}$ with $w = xy$.  Then $x\alpha_{s} > 0$,
$x\alpha_{t} > 0$.  We know from the discussion
of dihedral roots what the coefficients of
$$ y\alpha_{t} = p\alpha_{s} + q\alpha_{t} $$
are like, but then 
$$ w\alpha_{t} = p x\alpha_{s} + q x\alpha_{t} \; . $$
We can apply an induction hyothesis to $x\alpha_{s}$ and $x\alpha_{t}$
since $\ell(x) < \ell(w)$.\qed

\satz{Corollary}{Z}{The smallest positive root coefficient is
$1$.  Any root coefficient lying between $1$ and $2$ equals
$2\cos(\pi/m_{s,t})$ for some $s$, $t$.}
\endsatz

\proof/.  If
$$ { \sin (\pi n/m) \over \sin(\pi/m) } \quad (1 \le n \le m/2) $$
lies between $1$ and $2$ then $n=2$, and the least 
possible value is $\sqrt{2}$.
Any root coefficient will be a sum of terms
$$ C = c_{1} \dots c_{k} $$
where $k \ge 1$ and each $c_{i}$ will be either a positive
integer or one of the constants in the Proposition,
with each $c_{i} > 1$.  The smallest possible value
for one of these constants is $2 \cos (\pi/4) = \sqrt{2}$,
so if this product is less than $2$ then $k=1$,
and is an integer only if the product itself, of course,
is $1$.\qed

\satz{Corollary}{AA}{The smallest possible coefficient of a positive root,
other than $1$, is $\sqrt{2}$.}
\endsatz


\Section{Composition and decomposition}
One of the important results of
Brink's thesis is 
that the apparent difficulties concerning mixed cyclotomy 
in the determination of the minimal roots are
in fact spurious.
In order to see that it is true,
we have to look carefully at how roots are calculated.

A root $\sum \lambda_{\alpha} \alpha$ may be considered
as a function on the Coxeter graph:
$s \mapsto \lambda_{\alpha_{s}}$.
The {\bold support} of a root $\lambda = \sum \lambda_{\alpha} \alpha$
is its support as a function---the set of
$s$ where $\lambda_{\alpha_{s}} \ne 0$.  
Brink's nice idea is to track the computation of roots by 
looking at their support.  All roots are
constructed by applying some $w$ to
a basic root.   The support of a root $\lambda$
is extended by applying a reflection $s$ which is
not already in its support.  If $s$ is
not connected to that support in the Coxeter graph,
then $s\lambda = \lambda$.   This process gives a new root
only when $s$ is linked to the support of $\lambda$.  Therefore:

\satz{Proposition}{O}{The support of any root is connected.}
\endsatz

Suppose we apply a reflection $s$ to a positive root $\lambda$
whose support does not contain $s$.  Let $\alpha = \alpha_{s}$.
Only the coefficient $\lambda_{\alpha}$ changes, and to
$$
-\lambda_{\alpha} + 
\sum_{\lambda_{\beta} \ne 0, \, \alpha \sim \beta} [ -2\, \alpha \Dot \beta ] \lambda_{\beta} \; . $$
If $\lambda$ is a minimal root with $\lambda_{\alpha_{s}} = 0$, then $s\lambda$ will also
be a minimal root if and only if
$$ \sum_{\lambda_{\beta} \ne 0, \, \beta \sim \alpha} 
	[ - 2 \, \alpha \Dot \beta ] \lambda_{\beta} \le 2 \; . $$
But  $\lambda_{\beta} \ge 1$ and also
$- 2 \, \alpha \Dot \beta \ge 1$.  The condition for minimality
cannot hold here unless there is just one $\beta = \alpha_{t}$ in the sum.
Furthermore, in that case the above condition amounts to the condition
$$ [ - 2 \, \alpha \Dot \beta ] \lambda_{\beta} \le 2 \; . $$
The first term is either $1$ or $\ge \sqrt{2}$, as is the coefficient
$\lambda_{\beta}$.  All in all:

\satz{Proposition} {P}{Suppose $\lambda$ to be
a minimal root, $s$ an element of $S$ not
in the support $\Theta$ of $\lambda$.  If $s\lambda$ is
again minimal then there must be just a single
link from $s$ to $\Theta$, say to $t$, and one of these holds:
\smallbreak
\beginitems{3em}
\item{(a)} $m_{s,t} = 3$ and $\lambda_{\alpha_{t}} < 2$;
\item{(b)} $3 < m_{s,t} < \infty$ and $\lambda_{\alpha_{t}} = 1$.
\enditems
}
\endsatz

We can summarize this in a mnemonic diagram:
$$ \eqalign {
\bullet \quad & x \dash \, 0 \reflectsto x \dash x \quad (x < 2) \cr
\bullet \quad & 1 \Dash \, 0 \reflectsto 1 \Dash \, c \quad \big( c = 2\cos(\pi/m) \big) 
} $$

\satz{Corollary}{BB}{The support of a minimal
root is a tree with containing no links of infinite degree.}
\endsatz

One of Brink's good ideas is to
focus on the coefficients of a minimal
root which are equal to the minimum possible value $1$.  
Suppose $\lambda$ to be a
minimal root with $\lambda_{\alpha} = 1$.  
We know that $\lambda$ has been built up
from a chain
of minimal root predecessors 
$$ \lambda_{1} \prec \dots \prec \lambda_{n} = \lambda \; . $$
Some one of these, say $\mu = \lambda_{k}$, will be the first
to have a non-zero coefficient
at $\alpha$, and that coefficient must be $1$,
since it can only grow.
Unless $\mu = \lambda_{1} = \alpha$, 
the root $\mu$ must have
been constructed by extension
from a root $\lambda_{k-1}$ whose support did not contain
$s$, and the earlier Proposition tells
us that $\mu = s \lambda_{k-1}$,
where $s$ is linked to the support 
of $\lambda_{k-1}$ at a single node $t$ by
a simple link in the Coxeter graph,
and $\lambda_{k-1, \alpha_{t}} = \lambda_{k, \alpha_{t}} = 1$.
We have
$\lambda = w\mu$ where $w$ is
a product of generators other than $s$.

Since the support $\Theta$ of $\lambda$ is
a tree, the complement of $\alpha$ in $\Theta$ is the disjoint
union of exactly as many components $\Theta^{\bullet}_{i}$ as 
there are edges in the support
of $\lambda$ linked to $\alpha$, say $\ell$.  Reflections corresponding to nodes in
different components will commute.   Hence
we can write $\lambda = w_{1} \dots w_{\ell} \mu$ where
each $w_{i}$ is in $W_{\Theta^{\bullet}_{i}}$.
All the $w_{i}$ commute, so we can write 
the product in any order.
Each of the $\lambda_{i} = w_{i}\alpha$ is itself a root with support
in $\Theta_{i} = \Theta^{\bullet}_{i} \cup \{ \alpha \}$, minimal
because it will be a predecessor of $\lambda$ in the root graph.
We have therefore proven:

\satz{Proposition}{Q}{Let $\lambda$ be a minimal root with $\lambda_{\alpha} = 1$.
Let $\Theta$ be the support of $\lambda$,  the $\Theta^{\bullet}_{i}$ the connected
components of\/ $\Theta \backslash \{ \alpha \}$,
$\Theta_{i} = \Theta^{\bullet}_{i} \cup \{ \alpha \}$.  
The restriction
$\lambda_{i}$ of $\lambda$ to $\Theta_{i}$ is also
a minimal root.}
\endsatz

This is called the {\bold decomposition} of $\lambda$.
Conversely:

\satz{Proposition}{R}{Let $\Theta$ be a subset of $S$, $\alpha$
an element of\/ $\Theta$, the $\Theta^{\bullet}_{i}$ the components of
$\Theta \backslash \{ \alpha \}$, $\Theta_{i} = \Theta^{\bullet}_{i}\cup \{ \alpha\}$.
If for each $i$ the minimal
root $\lambda_{i}$ has support on $\Theta_{i}$ with coefficient
$1$ at $\alpha$, then the vector $\lambda$ with support $\Theta$
whose restriction to each $\Theta_{i}$ agrees with $\lambda_{i}$
is a minimal root.}
\endsatz

This is called {\bold composition}. 

\proof/.  More difficult than the previous result.
It follows from this result of Brink's:

\satz{Lemma}{F}{If $\lambda $ is any positive root with coefficient $1$
at $\alpha$, then $\lambda = w\alpha$, where $w$ is a product
of $s$ in $S$ not equal to $s_{\alpha}$.}
\endsatz

\proof/.  Induction on the depth of $\lambda$.
Trivial if $\lambda = \alpha$.  It remains to be shown that
if $\delta(\lambda) > 1$ and $\lambda_{\alpha} = 1$ then we can find
$s$ such that $s\lambda \prec \lambda$ and $s\lambda$
also has coefficient $1$ at $\alpha$.   There exists
at any rate some $s$ with $s\lambda \prec \lambda$.  If this
does not have coefficient $1$ at
$\alpha$, then---as we have seen above---$s\alpha$
has coefficient $0$ at $\alpha$;  $\alpha$ has exactly one link
to the support of $s\lambda$, say to $\beta$; and the coefficient of
$s\lambda$ at $\beta$ is equal to $1$.  By induction,
$s\lambda = w\beta$ where $w$ is a product
of $t$ equal neither to $s$ nor $s_{\beta}$.  But then
$$ \lambda = s w \beta = ws \alpha \; . \qed $$

These results have a natural generalization.
The {\bold unit support} of $\lambda$ is the set
of $\alpha$ with $\lambda_{\alpha} = 1$.
Suppose $\lambda$ to be a minimal root,
$\Theta$ its support, 
$T \subset \Theta$ its unit support.
We know that $\Theta$ is a connected tree.
The complement $\Theta \bs T$ can be expressed as the disjoint union
of its connected components $\Theta^{\bullet}_{i}$.
For each of these, let $\Theta_{i}$ be its union
with elements of $T$ attached to it.  The set
$\Theta_{i}$ will then be connected, and no
point of $T$ will be in its interior.
The set $T$ will contain in addition
some pairs of nodes of the Coxeter graph
connected by links in it; let
these be counted as additional
sets $\Theta_{i}$ (but with empty interiors).
Finally, there is one exceptional case
where $\Theta$ has just one element,
which is $\Theta_{1}$.  

\satz{Lemma}{G}{The set $\Theta$ is the union
of the $\Theta_{i}$, which overlap
only in points of $T$.}
\endsatz

\proof/.  We may suppose
$\Theta$ to have more than one element in it.
Suppose given $t$ in $\Theta$.  If it is not in
$T$, it belongs to a unique $\Theta^{\bullet}_{i}$.
Otherwise, it will be in $T$.  If it is connected in the Coxeter diagram
to one of the $\Theta^{\bullet}_{i}$, it will be in $\Theta_{i}$.
Finally, we have the situation where the
only elements of $\Theta$ it is linked to are in $T$.
Each of the pairs forming such a link
will be one of the $\Theta_{i}$.\qed

If $T$ is any subset of $S$, a subset
$\Theta$ is called {\bold $\bold T$-connected} if none of its interior points
belong to $T$.  Each of the $\Theta_{i}$ defined
above are $T$-connected, and are called its {\bold $\bold T$-components}. 
 Following Brink,
I call a minimal root {\bold indecomposable}
if its unit support is
contained in the boundary
of its support.

{\bold Lemma.} (Brink) {\sl (Decomposition) The restriction of
a minimal root to any of the $T$-components of its support
is again a minimal root.}

\proof/.  This follows by induction from the earlier decomposition.\qed

\satz{Proposition}{S}{{\rm(Brink)(Composition)} Let $\Theta$ be any
connected, simply connected
subset of $S$, $T$ an arbitrary subset of $\Theta$,
$\Theta_{i}$ its $T$-components.
If $\lambda_{i}$ is for each $i$
a minimal root with support $\Theta_{i}$
and unit support $T \cap \Theta_{i}$,
then the vector $\lambda$ with support $\Theta$
whose restriction to $\Theta_{i}$ is $\lambda_{i}$ will
be a minimal root.}
\endsatz

This reduces the calculation of minimal
roots to the calculation of
indecomposable ones, together
with a bit of accounting.  Finding
the reflection table turns out
not to be much more difficult.

How are indecomposable roots found?
They are obtained by reflection from predecessors in the root
graph.  What can we say about this process?
If $\lambda$ is an indecomposable root and $\mu = s_{\alpha}\lambda \prec \lambda$,
what can we say about $\mu$?  Let $\Theta$ be the support of $\lambda$.
There are several possibilities.
(1) Suppose that $\alpha$ is a boundary point of $\Theta$, and that $\lambda_{\alpha} = 1$.
Then  $\lambda$ is obtained from $\mu$
by {\bold extension} from a unit node of $\mu$.  This must
be a boundary point of the support of $\mu$.  This can happen only if
$\mu$ has singleton support.
(2) Suppose that $\alpha$ is in the boundary of $\Theta$ but that
$\lambda_{\alpha} > 1$.  Then $\mu$ will also be indecomposable.
(3) Suppose $\alpha$ lies in the interior of $\Theta$
and $\mu_{\alpha} > 1$.
Then $\mu$ is also indecomposable.
(4) Finally, suppose $\alpha$ lies in the interior
of $\Theta$ but $\mu_{\alpha} = 1$.  The support
of $\mu$ is also $\Theta$, so $\mu$ cannot be indecomposable.  It will
the composite of several $\mu_{i}$ adjoined to each other
at $\alpha$.  Thus:

\satz{Proposition}{T}{An indecomposable root
is constructed through a sequence made up of three types of processes:
(1) {\bold extension}, where $\mu \prec \lambda = s\mu$,
$\mu$ is also indecomposable, 
and the support of $\lambda$ is one node larger than
the support of $\mu$;
(2) {\bold fusion}, where $\mu$ is the composite of
indecomposable roots $\mu_{i}$ all with
$\alpha$ in the unit support;
(3) {\bold promotion}, where $\mu$ is indecomposable
and has the same support as $\lambda$.
}\endsatz

We have already investigated extension.
In the next section we look at fusion.

\Section{Reflecting at junctions}
Suppose $\Theta$ to be a connected and simply connected subset of $S$,
$\alpha \in \Theta$, the $\Theta_{i}$ the connected components
of $\Theta\bs \{s\}$ joined with $\{ s \}$.  Suppose that $\lambda$ is a minimal
root with support $\Theta$ and $\lambda_{\alpha_{s}} = 1$.
Let $\lambda_{i}$ be
the restriction of\/ $\lambda$ to $\Theta_{i}$, and 
assume the $\lambda_{i}$ are indecomposable.
{\sl Under what circumstances is $s\lambda$
a minimal root?}

Each of the $\Theta_{i}$ contains
a unique node $t_{i}$ linked to $s$.
Let $\beta_{i} = \alpha_{t_{i}}$.

The inner product $\lambda \Dot \alpha$ is 
$$ 1 - \sum  \lambda_{\beta_{i}} \, | \alpha \Dot \beta_{i} | \; . $$
Recall that $\lambda_{\beta_{i}} \ge 1$, $| \alpha \Dot \beta_{i} | \ge 1/2$.
When the junction has more than $3$ roots
joining at it, this dot-product is $\le -1$,
and $s\lambda$ will not be minimal.  

\satz{Proposition} {U}{If $\Theta \bs \{ s \}$ has more than $3$ components,
$s\lambda = \oplus$.}
\endsatz

We may now assume there to be
$2$ or $3$ links at $\alpha$. 

\satz{Lemma}{H}{Suppose $c = 2 \cos (\pi/m_{s,t})$.  The 
links $c \dash 1$ and $1 \Dash 1$ cannot occur in any root.}
\endsatz

This introduces notation I hope
to be self-explanatory.
What I mean by this is that if\/ $\lambda$ is a root
then there do not exist neighbouring nodes
$s$ and $t$ in the Coxeter graph with
(1) $m_{s, t} = 3$, $\lambda_{\alpha_{s}} = 1$,
$1 < \lambda_{\alpha_{1}}< 2$;
or (2) $m_{s, t} > 3$, $\lambda_{\alpha_{s}} = 1$,
$\lambda_{\alpha_{1}} = 1$.

\proof/.
By induction.  The only way to get $c \dash 1$ is from $c \dash 1 \dash 1$.\qed

We now consider the possibilities for fusion.

{\sl $\bullet$ Two links.}

The three {\sl a priori} possibilities are
$$ \eqalign {
& \lambda_{x} \dash 1 \dash \lambda_{y} \cr
& \lambda_{x} \dash 1 \Dash \lambda_{y} \cr
& \lambda_{x} \Dash 1 \Dash \lambda_{y} \cr
} $$

\satz{Proposition}{I}{In the third case $s\lambda = \oplus$.}
\endsatz

\proof/. The inner product is
$1 - c_{x} \lambda_{x}/2 - c_{y} \lambda_{y}$.  
We know that each $c_{*}$ is at least $\sqrt{2}$,
and by the previous Lemma that $\lambda_{*} \ge \sqrt{2}$ as well.
Thus the dot product is $-1$ or less.\qed

In the first case, the condition is $1 - \lambda_{x}/2 - \lambda_{y}/2 > -2$,
or $\lambda_{x} + \lambda_{y} < 4$.  The case $\lambda_{*} = c$ is excluded,
so either $\lambda_{*} = 1$ or $\lambda_{*} \ge 2$.  It is not
possible for both $\lambda_{*}$ to be $\ge 2$,
so one at least must be $1$, say $\lambda_{y} = 1$.

We can do better:
the inner product must be $-\cos(\pi k/m)$ for $k \le m/2$.  This gives
us: 

$m=3$.  The inner product is $0$,
with
$$ \bullet\;\; \quad 1 \dash \,1\, \dash \,1 \reflectsto  1 \dash \,1\, \dash \,1 $$
or $-1/2$ with
$$ \bullet\;\; \quad 2 \dash \,1\, \dash \,1 \reflectsto 2 \dash 2 \dash \,1 \; . $$

$m=4$.  We must have
$1/2 - \lambda_{x}/2 = - \sqrt{2}/2$ or $\lambda = 1 + \sqrt{2}/2$.
Impossible.

$m=5$.  $1 - \lambda_{x} = -c$ or $1-c$, $\lambda_{x} = 1 + c$
or $\lambda = c$.  The second is impossible, alas.  So the only possibility is
$$ \bullet\;\; \quad (1+c) \dash \,1\, \dash 1 \reflectsto (1+c) \dash \,(1+c)\, \dash 1 \; . $$

$m \ge 6$.  Nothing possible.

Now the remaining case
$$ \lambda_{x} \dash \,1\, \Dash \lambda_{y} \; , $$
Here 
$$ 1 - \lambda_{x}/2 - (c/2) \lambda_{y} > -1, \quad
\lambda_{x} + c \lambda_{y} < 4 \; . $$

$m=4$.  Only possibility is
$$ \bullet\;\; \quad 1 \dash 1 \Dash \sqrt{2} \reflectsto 1 \dash \,2\, \Dash \sqrt{2} \; . $$

$m=5$.  Only possibility is
$$ \bullet\;\; \quad 1 \dash \,1 \Dash \,c\, \reflectsto 1 \dash \,(1+c)\, \Dash \,c \; . $$

Note that in all these cases we have a join of two
primitive roots, one of which is just $1 \dash 1$.
So it is trivial to reconstruct the reflection even
when this is all embedded in
a generic root.  

{\sl $\bullet$ Three links.}

Assume the reflection is a primitive root.
The conditions are $\lambda_{x} + \lambda_{y} + \lambda_{z} < 4$ if
all links in are simple, or
$\lambda_{x} + \lambda_{y} + c_{m} \lambda_{z} < 4$ if one has degree $m > 3$.
More of degree $m > 3$ cannot occur by simple calculation.  If\/ $m > 3$ then
$c \ge \sqrt{2}$ so the second sum is at least $4$, and cannot occur.

As for the first, we cannot have any $\lambda_{*} \ge 2$,
so each can only be $1$ or $c$.  But the link
$c \dash 1$ cannot occur.

\satz{Proposition} {V}{The only triple junction
producing a minimal root be reflection
at the junction is the case that all coefficients are $1$.  
It reflects to a central node with value $2$.}
\endsatz
\medbreak

\overfullrule=0pt

\captionedFigure{3truein}{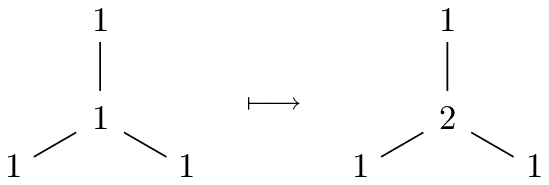}{The only triple junction
giving rise to fusion.}

\medbreak

Looking back over all ways indecomposable minimal
roots are constructed, we deduce:

\satz{Proposition}{CC}{No indecomposable minimal root
has more than one link
of degree $m > 3$ in its support.}\endsatz

This means that mixed cyclotomy does not occur.
Brink's thesis at this point goes on
to build a rather explicit description of all
possible indecomposable roots whose support
is a given Coxeter graph.  It is probably a poor
idea to try to use these lists
in a computer program, since it is fairly simple
to use basic properties of indecomposable roots
to build these lists automatically.
There is one consequence of what we have seen so
far that is simple to formulate.  It turns out that
any root coefficient may be described in no more
than three integers: each one is either of
the form $a + bc_{m}$ with $c_{m} = 2 \cos(\pi/m)$ or
larger than $2$.  With a bit of fiddling to handle
the last case carefully, the three integers
$(a, b, m)$ are sufficient.  This is a drastic version
of the principle of `no mixed cyclotomy'.  Implementing this idea
in an explicit program is straightforward enough,
but the program to be sketched in the next section,
which uses composition and decomposition in a more sophisticated fashion,
seems to be more efficient.

There is one final point to take into account.
When an indecomposable root $\lambda$
has degree $m > 6$,
it will have a unique link in its support of that 
degree.  The discussion of extension
above tells us that this link
is in the support of the
dihedral root first encountered in constructing $\lambda$.
Suppose this dihedral root has support
at $s$ and $t$, and that it started with $s$.
The discussion of dihedral roots tells us
that the coefficient of $\alpha_{s}$
in $\lambda$ will be either $1$ or greater than $2$.
In the first case, $s$ has to be on the boundary
of the support of $\lambda$ because of indecomposability.
In the second, we know that it is on the boundary
because no node with coefficient $\ge 2$
can be extended.   In all cases,
therefore, {$s$ is on the boundary of the support
of $\lambda$}.  I leave it as an exercise to show further
that {\sl in the second case the coefficient of
$\alpha_{t}$ is equal to $c_{m} = 2\cos(\pi/m_{s,t})$, and in both cases
all other coefficients are of the form $nc_{m}$ where $n$ is
a positive integer}.  (This, too, is an observation of Brink's
thesis.)


\Section{The program}
I am going to sketch here a program
that lists all minimal
roots and constructs the reflection table.
Relying strongly on Brink's thesis,
it works with two kinds of root,
{\bold indecomposable} and {\bold composite}.
The associated data structures are quite different.
Both have indices, descent sets,
and a reflection list.  But an indecomposable
root has in addition an actual array of coefficients.
These coefficients will all lie in
one of the real cyclotomic rings
generated by $\zeta^{\phantom{-1}}_{2m} + \zeta_{2m}^{-1}$,
where $m = m_{s,t}$ for some $s$, $t$.
I call these rings {\bold elementary}.
An indecomposable root is also assigned its degree $m$,
and when $m > 6$ its two special nodes, the ones
spanning the unique link of degree $m$, are specified as well.
A composite root has, instead,
a list of its indecomposable components.

New roots are found by applying
a reflection to a root already defined.
There are several basic ways in which to do
this:
(1) {\bold extension} of an indecomposable root to
produce a new indecomposable root with a larger support;
(2) {\bold promotion} of an indecomposable root to
produce a new indecomposable root with the same support;
(3) {\bold fusion} at a junction of several
indecomposable roots to produce
a new indecomposable root;
(4) {\bold composition} of an arbitrary root 
with a dihedral root attached at one
end to its unit support;
(5) {\bold replacement} of one or more of the components
of a composite root by a single indecomposable
root.

It is only the first three that require actual
arithmetic, and always in one of the
elementary real cyclotomic rings.  
Composition is an extremely simple operation,
while replacement amounts to
reducing a composite reflection to one of
the first three.  That this can be done easily
seems almost an accident depending
on the limited number of ways fusion can take place.
I'll not give details of this below, but I'll
sketch an example here.  Suppose we want to calculate
$s_{3}\lambda$ where 
$$ \lambda\Colon 1\, \dash \, 1 \, \dash \,1\, \Dash \sqrt{2} \dash \sqrt{2} $$
and $s = s_{3}$ is the third node from the left.  The root $\lambda$ is
the composite of three indecomposable roots
$$ 1\, \dash \, 1, \quad 1 \, \dash \,1, \quad 1\, \Dash \sqrt{2} \dash \sqrt{2}  $$
There are two components in the {\bold star} of $s$, 
the components of $\lambda$ containing it, and they
do not exhaust $\lambda$.  Instead, they make up
the smaller root
$$ \mu\Colon 1 \, \dash \,1\, \Dash \sqrt{2} \dash \sqrt{2}  $$
whose reflections we shall have already calculated
when $\lambda$ is taken off the queue.  It is easy enough to locate $\mu$,
given the reflection tables already made,
since it is equal to $s_{2}$ applied to the indecomposable
root
$$ \nu\Colon 1\, \Dash \sqrt{2} \dash \sqrt{2}  $$
We look up $s_{2}\nu$ to get $\mu$,
then look up $s\mu$, which is the indecomposable root
$$ 1 \, \dash \,2\, \Dash \sqrt{2} \dash \sqrt{2}  \; . $$
Finally, we replace the two components
of $\mu$ by the single component $s\mu$ in
the list of components of $\lambda$.  All this
can be carried out, with minor modifications,
for each of the cases of fusion.  In the pseudo-program
below the details will be grossly telescoped.

Here is the program:

First of all we define the basic roots, those in $\Delta$.

Then we define specially the dihedral roots,
those whose support consists of exactly two
nodes of the Coxeter graphs.  When we are through
with these we have finished the basic roots.
We put the dihedral roots of depth $2$ in the queue,
and after we have defined all of those we put
those of depth $3$ in the queue.  In this way,
we preserve the essential property of the queue
that roots of depth $d$ are put on before those of depth $d+1$.

As for dihedral roots of depth $4$ or more, we can finish them
as soon as they are constructed.  Both of their
coefficients are $2$ or more,
which means they cannot be extended, so the only
reflections higher in the root graph
are to dihedral roots with the same support.
In other words, we never put them
in the queue.  Incidentally, these
occur only for links with $m \ge 6$.
One consequence is that in calculating reflections
of roots taken off the queue, we have only to
do arithmetic with algebraic integers
of the form $a + b c_{m}$,
where $m$ is the degree of the root, $c_{m} = 2 cos(\pi/m)$,
$a$ and $b$ are integers.

At this point, all roots of depth $1$ have been finished, all
dihedral roots have been defined, and some have been finished.
The queue contains only dihedral roots.  While the queue is not empty,
we remove roots from it to be finished.

Let $\lambda$ be a root removed from the queue.
All roots of less depth have been finished already,
and in particular all descents $s\lambda \prec \lambda$ have
been assigned.  We run through all the $s$ in $S$,
and for each one of them where $s\lambda$ has not been 
assigned, we have to calculate $s\lambda$.  Exactly what this involves
depends on a number of things---the basic idea is to use
composition and decomposition to avoid repeating actual arithmetic.

So now suppose we are looking at $s$ and $\lambda$ with $s\lambda$ not known.
Let $\alpha = \alpha_{s}$.
There are many cases to deal with.

\medbreak
\Item{1em}{The node $s$ lies in the support of\/ $\lambda$.}
\Item{2em}{It belongs to a unique component $\lambda_{i}$.}
\Item{3em}{That component is $\lambda$ itself, which is therefore indecomposable.
We must calculate the dot product $\lambda \Dot \alpha$ to see if
$s\lambda$ is minimal or not, and if it is define it.}
\Item{3em}{The component $\lambda_{i}$ is not all of\/ $\lambda$, which is
a composite root.  We have calculated $s\lambda_{i}$ already,
and $s\lambda$ is obtained by replacing
$\lambda_{i}$ with $s\lambda_{i}$.}
\Item{2em}{It belongs to several components.
It is therefore a junction, and we are looking at a possible fusion.
Let $\mu$ be the star of $s$, the composite of those components whose supports contain $s$.}
\Item{3em}{The root $\mu$ is the same as $\lambda$.
We have to calculate $s\lambda$ explicitly.
It may be $\lambda$, $\oplus$, or a new 
indecomposable root.}
\Item{3em}{The root $\mu$ makes up only a part of\/ $\lambda$.  We have calculated
$s\mu$, and obtain $s\lambda$ by replacing $\mu$ with $s\mu$.
Exactly what we do to find $\mu$ and $s\mu$ depends on 
which case of fusion we are dealing with,
as I mentioned in the example I looked at earlier.
This can be decided by local information around $s$.}
\Item{1em}{The node $s$ does not lie in the support of\/ $\lambda$.}
\Item{2em}{It is linked to more than one node in the support,
or by a link of infinite degree.  Here $s\lambda = \oplus$.}
\Item{2em}{It is linked by a single link
of finite degree.}
\Item{3em}{The link is to a unit node of\/ $\lambda$.
Then $s\lambda$ is a new composite.}
\Item{3em}{The link is not to a unit node.}
\Item{4em}{The root $\lambda$ is indecomposable.
We have to calculate explicitly.}
\Item{4em}{The root $\lambda$ is composite.
If\/ $\mu$ is the component it is linked to,
we have already calculated $s\mu$.
We replace $\mu$ by $s\mu$ in $\lambda$.}


\newcount\refno
\refno=0
\def\refitem{\medskip
\advance\refno by 1
{\bold \the\refno. }}
\def\jour#1{{\frenchspacing \sl #1}}
\def\pp#1{{\frenchspacing #1}}

\def\book#1{{\sl #1}}
\def\vol#1{{\bold #1}}

\references
\refitem
Brigitte Brink, `The set of dominance-minimal roots',
available as Report 94--43 from the School of mathematics
and Statistics at the University of Sydney:

{\tt http://www.maths.usyd.edu.au:8000/res/Algebra/Bri/dom-min-roots.html}

\refitem
Brigitte Brink, `The set of dominance-minimal roots',
\jour{Journal of Algebra} \vol{206} (1998), \pp{371--412}.
This publication, although it has the same title as the previous
item, is very different from it.   For purposes
of computation it is decidedly less interesting, since it
doesn't introduce the indecomposable minimal roots,
and emphasizes hand-lists instead.

\refitem
Brigitte Brink, `On centralizers of reflections in Coxeter groups', 
published in \jour{Bulletin of the London
Mathematical society} \vol{28} (1996), \pp{465--470}.
Also available as a preprint from 

{\tt http://www.maths.usyd.edu.au:8000/res/Algebra/Bri/}

\refitem
Brigitte Brink and Robert Howlett, `A finiteness property and an automatic structure for
Coxeter groups', \jour{Math. Ann.} \vol{296} (1993), \pp{179--190}.

\refitem
Bill Casselman, `Automata to perform basic calculations in
Coxeter groups', in \book{Representations of groups},
{\sl CMS Conference Proceedings 16}, A.M.S., 1994.

\refitem
Bill Casselman, `Computation in Coxeter groups---multiplication',
preprint, 2001.
To appear in the \jour{Electronic Journal of Combinatorics}.

\refitem
Bill Casselman, `Java source code for finding minimal roots', at

{\tt http://www.math.ubc.ca/people/faculty/cass/java/coxeter/roots/}

\refitem
Fokko Du Cloux, `Un algorithme de forme normale pour les groupes de Coxeter',
preprint, Centre de Math\'ematiques \`a l'\'Ecole Polytechnique,
1990.

\refitem
J. Tits, `Le probl\`eme des mots dans les groupes de Coxeter',
\jour{Symposia Math.} \vol{1} (1968), \pp{175--185}.

\refitem
E. B. Vinberg, `Discrete linear groups generated by reflections',
\jour{Math. USSR Izvestia} \vol{5} (1971), \pp{1083--1119}.

\bye